\documentclass{svproc}
\usepackage[numbers,sort&compress]{natbib}
\usepackage{cite}
\usepackage{booktabs}
\usepackage{array}
\usepackage{url}
\usepackage{amsmath}
\usepackage{empheq}
\usepackage{multirow}
\usepackage{cases}
\usepackage{xcolor}
\usepackage{indentfirst} 
\usepackage{float}

\begin{document}
\mainmatter              

\title{Multiobjective Logistics Optimization for Automated ATM Cash Replenishment Process}
\titlerunning{Multiobjective Optimization for ATM Cash Replenishment}  


\author{Bui Tien Thanh\inst{1} \and Dinh Van Tuan\inst{1} \and Tuan Anh Chi\inst{1} \and
Nguyen Van Dai\inst{1} \and Nguyen Tai Quang Dinh\inst{1} \and Nguyen Thu Thuy\inst{1} \and Nguyen Thi Xuan Hoa\inst{2,}\inst{3}*}

\authorrunning{Bui Tien Thanh et al.} 
%
\tocauthor{Bui Tien Thanh, Dinh Van Tuan, Tuan Anh Chi, Nguyen Van Dai, Nguyen Tai Quang Dinh, Nguyen Thu Thuy, Nguyen Thi Xuan Hoa}
%


\institute{School of Applied Mathematics and Informatics, \and
School of Economics and Management, \inst{3} BK Fintech, \\ 
Hanoi University of Science and Technology, Hanoi, Vietnam\\
*Corresponding author:
\email{hoa.nguyenthixuan@hust.edu.vn}
}

\maketitle 

\begin{abstract}
In the digital transformation era, integrating digital technology into every aspect of banking operations improves process automation, cost efficiency, and service level improvement. Although logistics for Automated Teller Machine (ATM) cash is a crucial task that impacts operating costs and consumer satisfaction, there has been little effort to enhance it. Specifically, in Vietnam, with a market of more than 20,000 ATMs nationally, research and technological solutions that can resolve this issue remain scarce. In this paper, we generalized the vehicle routing problem for ATM cash replenishment, suggested a mathematical model, and then offered a tool to evaluate various situations. When being evaluated on the simulated dataset, our proposed model and method produced encouraging results with the benefits of cutting ATM cash operating costs.
\end{abstract}
\section{Introduction}
 ATM is a physical interaction point and can be compared to an extension arm of banks to customers in many areas. Therefore, storing and replenishing cash for ATMs affect customer satisfaction with banking services. Besides, there are also costs associated with personnel, transportation and interest rate on money deposited at ATMs - this is a relatively high cost since this rate can be as high as $10\%$ \citep{link:WorldBank}. From the supply chain management perspective, this is a classic trade-off between transportation and inventory costs. Many cost-optimizing models have been devised to solve this problem, especially in manufacturing industries \citep{peres2017optimization}. For the banking sector, this can be considered a money supply chain problem and a specific commodity that must adhere to the State's regulations and risk management in logistics. According to the International Monetary Fund (IMF), Vietnam has 20,404 ATMs, with 2,644 ($\approx$13\%) belonging to Vietnam Bank for Agriculture and Rural Development (Agribank). With an extensive network covering many regions and continuous cash replenishment, banks can improve the logistics of the ATM cash replenishment process, which will result in cost and administration savings.

The ATM cash replenishment problem aims to provide a plan to meet the demand for ATM withdrawal and optimize the costs related to money transportation. 
Many mathematical models 
\citep{ekinci2015optimization,Kiyaei2021deepqnetwork} been proposed to solve the transportation optimization problem, also known as the Vehicle Routing Problem (VRP). The problem is commonly expressed as a single objective optimization issue, which involves reducing the total transportation expenses of vehicles. This cost includes fixed costs and costs proportional to the total distance travelled. We also consider financial costs an important factor in our objective function. To address the VRP, we formulate the problem as a Multiobjective Optimization Problem (MOP). MOP involves multiple objectives that may conflict with one another, 
making it impossible to achieve an optimal solution with respect to a single objective. Thus, the aim is to find solutions that balance different objective functions best. Some methods to solve MOP are described in 
\citep{thang2016solving,thang2020monotonic}. In Vietnam, there are few research and technological solutions for this problem. Especially when the Vietnamese banking sector is in the digital transformation, many internal business processes are being digitized but not integrated with optimal technologies. Therefore, an easy-to-apply, close-to-logistics model in the Vietnamese market to improve ATM operations can benefit domestic businesses.

In this study, we generalize the replenishment cash process for ATMs based on a survey in the domestic banking market to better define the problem and the desired goal. Then, we model the problem, propose a tool to solve the problem and perform tests on the simulated dataset to evaluate the results. We expect that our recommendations would benefit the banking industry by: (i) Reducing interest costs on money deposited for a long time at ATMs; (ii) Saving transportation costs and time since distances and trips are optimized; (iii) Minimizing planning time and errors; (iv) Utilizing a system and related toolkits to standardize the procedure and reduce the personnel-dependent risk; (v) Providing a long-term transportation plan and an overview to assist the management department in making decisions regarding cash flow coordination.

The rest of this paper is organized as follows. Section 2 presents the problem and business requirements. Section 3 models the problem. Our proposed solution is presented in section 4. Section 5 tests and evaluates the results. Finally, section 6 concludes our main contributions.
\section{Research problem}

\begin{figure}[htbp]
    \centering
    \includegraphics[width = 0.9\textwidth]{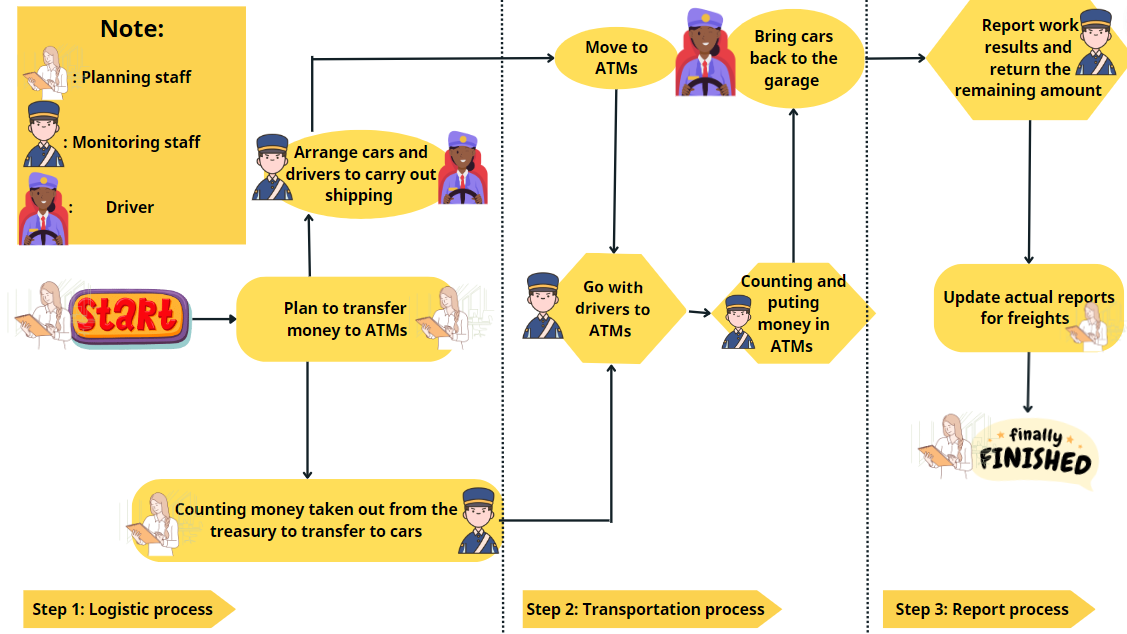}
    \caption{ATM cash replenishment process}

    \label{Fig1}
\end{figure} 

Fig. \ref{Fig1} shows a basic ATM replenishment process based on our survey \citep{serengil2019atm}.
As a first step, staff members plan to top-up the ATMs, i.e. determine which ATMs will be recharged, when to top up, and the amount of money. Once the plan is done, the bank's staff will coordinate with the shipping department or a third party to transfer the money to the ATMs. Fig. \ref{Fig1} only records those directly and full-time involved in replenishment; In fact, other management and support departments are also included, but less often. Generally, the cash replenishment process for ATMs has some common regulations, such as:
\begin{itemize}
    \item Money transportation needs to be done during office hours.
    \item An ATM is usually provided cash by a permanent warehouse.
    \item If a bank has many cash depots, usually the depots will be in charge of separate areas to meet the criterion of the shortest distance. This is also convenient for grouping personnel and managing cash flow by region.
\end{itemize}

And in the ATM cash replenishment process, the planning procedure plays an important role because it affects the incurred costs of the operation, including two main categories:
\begin{itemize}
    \item Firstly, transportation costs: are the costs related to renting a specialized vehicle to transport money, hiring a driver, salary and allowances paid to the staff planning, checking money and supervising the transport vehicle.    
    \item Secondly, financial cost: includes loan interest to keep money in ATMs to serve customers' withdrawal needs.
\end{itemize}

However, according to our survey in a segment of domestic banks in Vietnam, this process can be improved as shown in Table \ref{tab:Table2}.

\begin{table}
    \centering
    \caption{Planning for the ATM cash replenishment process}
    \vspace{-0.2cm}
    \label{tab:Table2}
    \begin{tabular}{ | m{6cm}| m{6cm} | } \hline
        \textbf{\qquad Current planning process} & \textbf{\qquad Things can be improved} \\ \hline
        Transportation planning is usually made in a short time, specifically 1-3 days, which means that the staff usually plans one day in advance and takes 2-3 days later to combine transportation. & Create a longer-term transportation planning (7-14 days) to find a better mix and provide an overall vision of the cash plan for management and other departments. \\ \hline
        The amount of money deposited for each ATM normally follows a fixed rate for a period of time (monthly or quarterly).& Consider forecasting actual withdrawals at ATMs by day to make a better reasonable deposit. \\ \hline
        Plans are often made using Microsoft Excel with manual calculations and staff experience.&-Use proven mathematical theories and models to create cost-optimized plans.
        \newline
        -Use tools and software to calculate more accurately, faster, to save time and effort. \\ \hline
    \end{tabular}
\end{table}
Withdrawal demand forecasting data will be used to calculate the plan to pour money into ATMs. Currently, in Vietnam, depending on the location, the ATM will be supplemented with a certain amount of money when the amount is below a safe threshold. This amount is fixed monthly or quarterly and assessed for changes the following month or quarter. Although this approach is simple, it has many disadvantages. As withdrawal demand can fluctuate over a short period of time, setting a deposit limit at an ATM can result in large interest costs on deposits when actual withdrawal needs differ significantly from the limit. On the other hand, if the demand for withdrawals increases sharply over time, the bank must replenish more times. We see an opportunity to combine transportation routes when planning flexibly according to demand fluctuations rather than a fixed amount of money. 

\section{Mathematical Model}
\subsection{Problem Statement}
\begin{figure}[htbp]
    \centering
    \includegraphics[width = 0.95\textwidth]{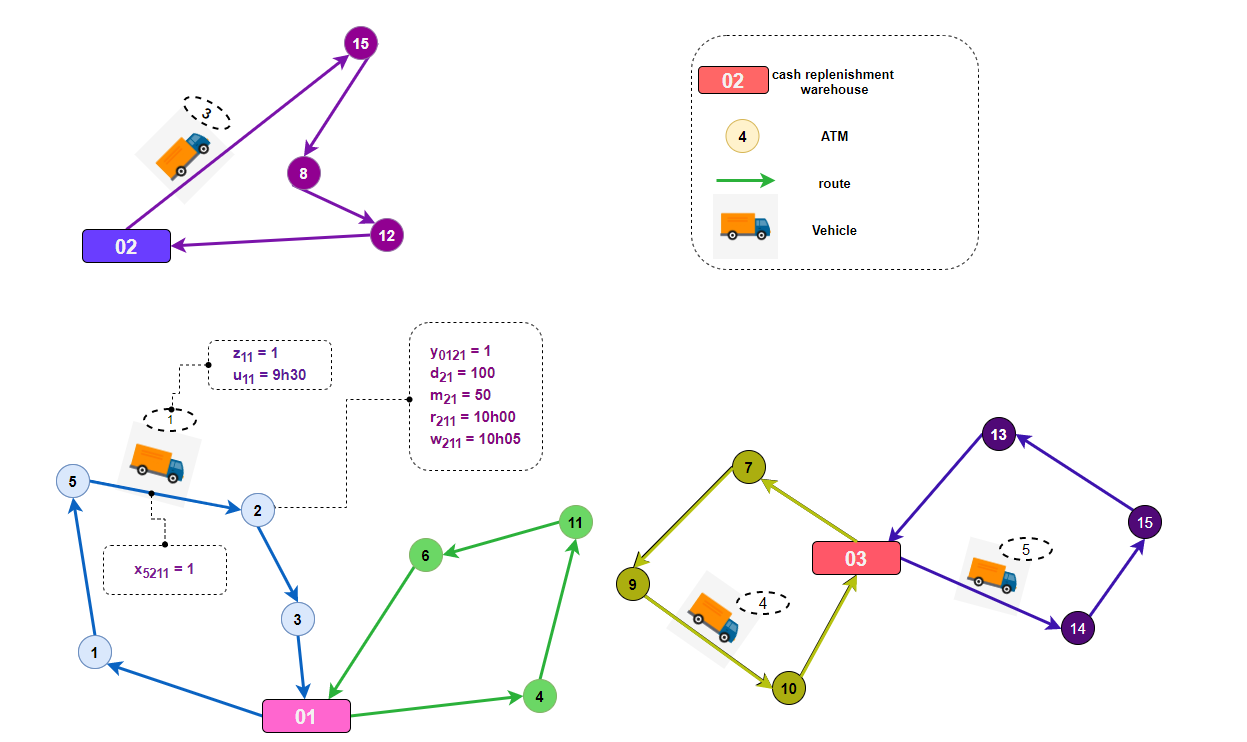}
    \caption{A scenario with 3 depots and 16 ATMs}
    \label{fig:figure_label}
\end{figure}

The problem is defined on a complete directed graph $\mathcal{G} = (\mathcal{V}, \mathcal{E})$, where $\mathcal{V} = \{1,2,\dots, A,01,02,\dots, 0D\}$ is the set of vertices of the graph and $\mathcal{E} = \{(i,j)| i,j\in\mathcal{V}, i\neq j\}$ is the set of edges of the graph, with $D$ and $A$ are positive integers. $\mathcal{A} = \{1,2,\dots, A\}$ is the set of ATMs and $\mathcal{D} = \{01,02,\dots, 0D\}$ is the set of cash depots. $\mathcal{H} = \{1,2,\dots, H\}$ is the set of vehicles. $\mathcal{T} = \{1,2,\dots, p\}$ is the set of transportation periods. Fig. \ref{fig:figure_label} illustrates an example of 3 depots and 16 ATMs. The ATMs $1,2,3,4,5,6,11$ are replenished from depot $01$. The ATMs $8,12,16$ are replenished from depot $02$. The ATMs $7,9,10,13,14,15$ are replenished from depot $03$. 
The optimal route for vehicle number 1 is $01 \to 1 \to 5 \to 2 \to 3 \to 01$. The optimal route for vehicle number 4 is $03 \to 7 \to 9 \to 10 \to 03$. Vehicle number 1 is used in period 1 $(z_{11}=1)$ and leaves the depot at 9h30 $(u_{11} = 9$h$30)$. ATM number 2 was replenished by depot 01 on period 1 $(y_{0121} = 1)$. The road (5,2) has vehicle number 1 passing on period 1 $(x_{5211} = 1)$. The amount to be replenished for ATM number 2 on period 1 is 100 million VND $(d_{21} = 100)$. The amount withdrawn at ATM number 2 on period 1 is 50 million VND $(m_{21} = 50)$. Vehicle number 1 arrives at ATM number 2 on period 1 at 10h00 $(r_{211} = 10$h$00)$. Vehicle number 1 started replenishing money at ATM number 2 on period 1 at 10h05 $(w_{211} = 10$h$05)$.

\subsection{Constraints} 
The model must satisfy the following constraints. (1) Each vehicle departs from and returns to the same depot. (2) In each period, each ATM can only be served by one depot. (3) In each period, each ATM can only be served by one vehicle. (4) The vehicles are heterogeneous (have different capacities and operating costs). (5) The start time of cash replenishment at each ATM must be within the allowable range. (6) Vehicles can only operate within specified hours. (7) The actual time taken at each ATM must be considered (to perform cash counting, cash loading, machine shutdown, et cetera). (8) The amount of cash replenished at each ATM in each period must not exceed the capacity of each vehicle. (9) The amount of cash to be transferred to each ATM is known in advance. (10) Only the total amount of cash replenished at ATMs is considered without considering the details of the denomination. (11) The location and coordinates of each ATM are known in advance. (12) ATMs function normally without breaking down or needing repair. 

\subsection{Mathematical Model} 
In this article, several parameters will be used to help model the problem, such as: $t_{ijh}$: The time for vehicle $h$ to travel from node $i$ to node $j$; 
$C$:The maximum travel distance allowed for vehicle $h$ within the period $t$; 
$q_h$: The capacity of vehicle $h$; 
$IR$: The annual interest rate $(\%)$;  
$I_{0j}$: The initial amount of money at ATM $j$;  
$c_{ij}$: The distance between node $i$ and node $j$; 
$e_j$: The earliest time to start replenishing money at ATM $j$; 
$l_j$: The latest time to start replenishing money at ATM $j$; 
$s_j$: The time needed to replenish money at ATM $j$ (measured from the time the money is taken out of the vehicle); 
$e_0$: The earliest time for the vehicle to start from the depot; 
$l_0$: The latest time for the vehicle to start from the depot; 
$a_h$: The driving cost of a vehicle $h$ per unit distance; 
$V$: The number of vertices in the graph ($V=|\mathcal{V}|$); 
$t_{max}$: The maximum travel time allowed for vehicle $h$ within the period $t$; 
$p$: The maximum number of periods allowed for all; 
$m_{jt}$: The amount of money withdrawn from ATM $j$ during $t$.

In addition, the article uses some decision variables to increase the accuracy of the problems such as: 
$x_{ijht}$ = 1 if vehicle $h$ travels on road $(i,j)$ during period $t$; 0, otherwise; 
$z_{ht}$ = 1 if vehicle $h$ is used during period $t$; 0, otherwise; 
$y_{ijt}$ = 1 if ATM $j$ is replenished by depot $i$ during period $t$; 0, otherwise; 
$d_{jt}$ = The amount of money needed to replenish ATM $j$ in period $t$; 
$r_{jht}$ = The time when vehicle $h$ arrives at ATM $j$ during period $t$; 
$w_{jht}$ = The time when vehicle $h$ starts replenishing money at ATM $j$ during period $t$ (measured from the time the money is taken out of the vehicle); 
$u_{ht}$ = The time when vehicle $h$ starts from the depot during period $t$.

Based on the situation discussed earlier, we have formulated a constrained optimization problem with the objective function represented as follows:
\begin{equation*}
\mathrm{min} \; \boldsymbol{f} = \mathrm{Vmin}\;(f_1, f_2)
\end{equation*}
Where Vmin is the minimum vector and
\begin{numcases}{} 
f_1 = \displaystyle{\sum_{i\in\mathcal{V}}\sum_{j\in\mathcal{V}|i\neq j}\sum_{h\in\mathcal{H}}\sum_{t\in\mathcal{T}}a_hc_{ij}x_{ijht}}\hspace{0.2cm} \label{obj_f1}\\
f_2 = \displaystyle{\sum_{j\in\mathcal{A}}\frac{IR}{365}\Big(p(I_{0j} + d_{j1}- m_{j1}) + (p-1)(d_{j2}- m_{j2}) + \cdots } \nonumber \\
\hspace{3.5cm} \begin{aligned}
&\quad +2(d_{jp-1}- m_{jp-1}) + (d_{jp}- m_{jp})\Big) \label{obj_f2}
\end{aligned}
\end{numcases}
s.t.
\begin{align}
    &\sum_{i\in\mathcal{V}}\sum_{h\in\mathcal{H}} {x_{ijht}}  = 1, \quad \forall j\in \mathcal{A}, t\in \mathcal{T} \label{con_in}\\
    &\sum_{i\in\mathcal{V}}\sum_{j\in\mathcal{A}} {d_{jt}} {x_{ijht}} \le q_h, \quad \forall h\in \mathcal{H}, t\in \mathcal{T} \label{con_demand}\\
    &\sum_{i\in\mathcal{D}}\sum_{j\in\mathcal{A}} {x_{ijht}} = \sum_{i\in\mathcal{D}}\sum_{j\in\mathcal{A}} {x_{jiht}} = {z_{ht}}, \label{con_use_vehicle}\quad \forall h\in \mathcal{H}, t\in \mathcal{T} \\
     &0 \le {r_{jht}} \le {w_{jht}}, \quad \forall j\in\mathcal{A},h\in\mathcal{H},t\in\mathcal{T}\label{con_arrival}\\
    & e_j \le {w_{jht}} \le l_j, \quad \forall j\in\mathcal{A},h\in\mathcal{H},t\in\mathcal{T}\label{con_start_serve} \\
 &\sum_{i\in \mathcal{S}}\sum_{j\in \mathcal{S}} {x_{ijht}} \le |\mathcal{S}| - 1, \quad \mathcal{S}\subseteq\mathcal{A}, 2 \le |\mathcal{S}| \le A, \forall h\in \mathcal{H}, t\in \mathcal{T} \label{con_subtour} \\
    &\sum_{i\in\mathcal{V}}\sum_{j\in\mathcal{V}|i\neq j}t_{ijh}{x_{ijht}} \le t_{max},  \quad \forall h\in \mathcal{H}, t\in \mathcal{T} \label{con_tmax} \\
    &\sum_{k\in\mathcal{V}}({x_{ikht}} + {x_{kjht}}) \le 1 + {y_{ijt}} 
   ,\quad\forall i\in\mathcal{D}, j\in\mathcal{A}, h\in\mathcal{H}, t\in\mathcal{T} \label{con_assign_atm}\\
 &\sum_{i\in\mathcal{D}} {y_{ijt}} \le 1, \quad \forall j\in \mathcal{A}, t\in\mathcal{T} \label{con_max_treasury}\\
    &\sum_{i\in\mathcal{V}}\sum_{j\in\mathcal{V}|i\neq j}\sum_{t\in\mathcal{T}}c_{ij}{x_{ijht}} \le C,\quad\forall h\in\mathcal{H} \label{con_cmax}
 \end{align}
\begin{align}
        &({w_{jht}}+ s_j + t_{ijh} - {r_{jht}}) {x_{ijht}} = 0,
    \quad\forall i\in\mathcal{V}, j\in\mathcal{A}, h\in\mathcal{H}, t\in\mathcal{T} \label{con_itoj} \\
    &{m_{jt}} \leq t(I_{0j}+{d_{j1}}-{m_{j1}})+(t-1)({d_{j2}}-{m_{j2}}) + \cdots \notag \\
&\quad\quad \quad\quad+ 2({d_{jt-1}}-{m_{jt-1}})+({d_{jt}}-{m_{jt}}) \label{con_withdraw} \\
     &\text{Let} \; T = \sum_{i\in\mathcal{V}}\sum_{i\in\mathcal{V}|i\neq j}t_{ijh}{x_{ijht}} + \sum_{i\in\mathcal{V}}\sum_{j\in\mathcal{A}|i\neq j}s_j{x_{ijht}}\notag \\
 & l_0 - T \le {u_{ht}}, \quad\forall h\in\mathcal{H}, t\in\mathcal{T} \label{con_depart} \\
 &T \le l_0 - e_0, \quad \forall h\in\mathcal{H}, t\in\mathcal{T} \label{con_duration} \\
 &{x_{ijht}}, {y_{ijt}}, {z_{ht}} \in \{0,1\}, \quad\forall i,j\in \mathcal{V},h\in\mathcal{H}, t\in\mathcal{T} \label{con_binary}\\
 &0 \le {d_{jt}} \quad \forall j\in\mathcal{A},t\in\mathcal{T}\label{con_non_negative}
\end{align}
The objective function (\ref{obj_f1}) minimizes transportation costs, while (\ref{obj_f2}) minimizes financial costs. Constraint (\ref{con_in}) ensures that each ATM can be passed by only one vehicle per period. Constraint (\ref{con_demand}) ensures that the amount of money transferred to each ATM per period does not exceed the capacity of each vehicle. Constraint (\ref{con_use_vehicle}) ensures that a vehicle must depart from its corresponding depot when used. Constraint (\ref{con_arrival}) ensures that the time a vehicle starts unloading money at ATM $j$ in period $t$ is no earlier than its arrival time. Constraint (\ref{con_start_serve}) ensures that the time a vehicle starts unloading money at ATM $j$ in period $t$ is within the working hours of ATM $j$. Constraint (\ref{con_subtour}) eliminates sub-tours. Constraint (\ref{con_tmax}) ensures that the travel time of each vehicle in each period is within the threshold. Constraint (\ref{con_assign_atm}) allows the replenishment of an ATM by a depot only if there is a route from the depot to the ATM. Constraint (\ref{con_max_treasury}) ensures that each ATM is replenished by at most one depot per period. Constraint (\ref{con_cmax}) ensures that the distance traveled by each vehicle in each period does not exceed the threshold. Constraint (\ref{con_itoj}) relates to the order of passing through vertices: if vehicle $h$ travels directly from vertex $i$ to $j$, the arrival time ${r_{jht}}$ at vertex $j$ must equal $({w_{iht}}+ s_i + t_{ijh})$. Constraint (\ref{con_withdraw}) ensures that the amount of money withdrawn at ATM $j$ in period $t$ does not exceed the amount currently available at ATM $j$. Constraints (\ref{con_depart}) and (\ref{con_duration}) ensure that a vehicle's operating time is within the designed working hours $[e_0, l_0]$. Constraint (\ref{con_binary}) indicates that decision variables only take the values 0 or 1. Constraint (\ref{con_non_negative}) ensures that variables are non-negative.

Pareto-based multiobjective learning algorithms can generate multiple Pareto-optimal solutions, allowing users to gain insights about the problem and make better decisions. To select the best Pareto-optimal solution from the Pareto front, sub-criteria might be employed and optimized accordingly. The mathematical model associated with this problem is the optimization problem on the efficient set, which has been studied in prior works \citep{thang2015outcome}.
In future research, we plan to build on these previous works to further explore and enhance the effectiveness of the optimization problem on the efficient set.

\section{Methodology}
In our survey, VRP is a long-standing problem with many studies conducted. These solutions are usually divided into two categories: the group of exact solution methods, typically the branching method with the mixed integer programming problem model \citep{golden2023evolution, aggarwal2019mixed}, or the group of metaheuristic methods such as evolutionary search locally algorithms, local search \citep{golden2008vehicle}.

When considering practical applications, transportation problems will often be reduced to planning problems and use a specific algorithm software such as CPLEX, MATLAB or AMPL to solve. In this article, we will use open-source software called OR-Tools. The most significant difference between OR-Tools and other software is that it has prepackaged components and constraints of some familiar optimization problems with popular functions that users can use directly. This makes the model simpler and easier to visualize when setting constraints, thereby saving product research and development time.

Some problems that OR\-Tools supports are assignment, scheduling, vehicle routing, and network flow problems. This study uses the tool of vehicle routing problems. Main steps in problem-solving by OR\-Tools are described below:
    Initialize data for the problem: information about the distance matrix between the points, the operating time of the depots and ATMs, the amount of money to be transferred to the points, information about the used vehicles...  
    In fact, we will balance the financial costs (interest on deposits in the ATM) and the transportation costs (travel costs of the vehicle) when transporting money to the ATM. Therefore, we perform a simple test by dividing a large money amount into several small amounts and estimating the financial cost proportional to the number of scheduled days. 
    Combine financial and transportation costs to obtain a cost matrix. The problem will be optimized based on this aggregate cost.
    Build models and load data using available functions in OR-Tools.
    Set the limit of the problem by 2 conditions: for vehicles and for ATMs.
    Set the algorithm, maximum program running time and run it.

After doing those steps, the result obtained from OR-Tools is a transportation plan for a specified period, with details of which points each transportation will pass in a day, the order of its traveling, and how much will be transported.

\section{Testing and Evaluation}
To perform the test, we use the data from actual ATM locations of a bank in Hanoi (Vietnam). The distance matrix between ATMs is collected from Google Maps to get the closest distance to reality, which includes the possible route, prohibited roads, ... Operating hours, stopping time at a point, and specialized vehicles transporting money are estimated and taken parameters based on an actual business survey at a domestic bank. For cost data, we refer to staff salaries, petrol, vehicle rental prices, and average deposit interest rates in the market.

We calculate the travel cost proportional to the distance and the financial cost depending on the amount with the estimated backlog days. The key issue is to balance the cost of money deposited at the ATM with the costs incurred from transporting money. To evaluate the effectiveness, we consider 2 cases:

(i) Case 1: Set a fixed amount of money to be deposited for ATMs and only perform cash replenishment for ATMs when this amount is nearly depleted.

(ii) Case 2: Divide ATM deposits into smaller amounts based on the forecasted changes in ATM withdrawals over a period of time.

We experiment with 2 different samples, including from 2 cash depots in different areas and within a 7-day period. The key issue to be addressed is to balance financial costs with transportation costs and fixed vehicle costs. We use the method of splitting large amounts of money into smaller amounts, simply called splitting orders. The idea is to determine the lower and upper bounds for small orders, and based on that calculate the number of split orders that can be made for the initial large order, and let the algorithm choose one of the calculated ways of splitting orders to optimize as much as possible. The experimental parameters and corresponding results are presented in Table \ref{tab:Table_5}.
\vspace{-0.5cm}
\begin{table}[!h]
    \begin{center} 
    \caption{ Test Results}
    \label{tab:Table_5}
    \begin{tabular}{ |m{3.4cm}|m{2.3cm}|m{1.8cm}|m{2.3cm}|m{1.8cm}| }\hline
       \multirow{2}{*}{ } & \multicolumn{2}{|c|}{Test 1} &\multicolumn{2}{|c|}{Test 2} \\ 
       \cline{2-5}
             & No single split & Single Split & No single split & Single split \\ \hline
        \multicolumn{5}{|l|}{Input} \\ \hline
        Number of warehouses & \multicolumn{4}{c|}{2 warehouses} \\ \hline
        Number of vehicles &  \multicolumn{4}{c|}{2 vehicles per warehouse} \\ \hline
        Number of days &  \multicolumn{4}{c|}{7 days} \\ \hline
        Approximately the amount offered for the entire period &  \multicolumn{4}{c|}{Ranges from 2.5 to 3.5 billion VND per ATM} \\ \hline
        Approximately adjusted money amount pouring each time & \multicolumn{4}{c|}{Ranging from 1.0 - 1.4 billion VND per ATM per provide} \\ \hline
        Bank interest rate &\multicolumn{4}{c|}{5\% per year } \\ \hline
        Number of ATMs & \multicolumn{2}{c|}{28} &\multicolumn{2}{c|}{58}\\ \hline 
        Total amount offered & \multicolumn{2}{c}{84.412.000.000 VND} & \multicolumn{2}{|c|}{176.999.000.000 VND}\\ \hline
        \multicolumn{5}{|l|}{Results when running data with OR-Tools} \\ \hline
        Total number of trips & 5 & 8 & 7 & 12 \\ \hline
        Total distance (km) & 326,8 & 621,9 & 520,4 & 976,9 \\ \hline
        Transportation cost (VND) & 5.707.054 & 16.181.443 & 10.878.837 & 31.078.885 \\ \hline
        Financial cost (VND) & 36.317.446 & 14.748.057 & 72.733.163 & 20.986.315 \\ \hline
        Total cost (VND) & 42.024.500 & 30.929.500 & 83.612.000 & 52.065.200 \\ \hline
        Cost improvement (\%) & \multicolumn{2}{c|}{26,4} &  \multicolumn{2}{c|}{37.7} \\ \hline      
    \end{tabular}
   \end{center}
\end{table}
Table \ref{tab:Table_5} shows that by splitting the number of replenishment according to the forecasting, the total cost, including transportation and financial costs, is improved by 26 - 38 \% with tested data. Although splitting single transportation into multiple ones with smaller amounts of money will increase the vehicle working time and the total distance and transportation costs, the cost of deposits at ATMs will decrease and offset the total costs. Accordingly, there is potential for using optimized models and tools to lower the bank's operating costs.

\section{Conclusion}

In this article, we generalized the logistics problem of cash replenishment for bank ATMs and presented a model and tool to support cash replenishment planning. The article may have some limitations, such as not considering actual costs that may arise in operation, such as indirect personnel costs, opportunity costs of investment cash flow, as well as not considering special cases: Planning for holidays (when demand suddenly increases), ATM breakdowns and other technical problems. However, we expect our recommendations to bring benefits in procedure automation, management, and operating costs for banks in Vietnam.

%

\renewcommand\bibname{References}
\bibliographystyle{splncs03}
\bibliography{references.bib}

\end{document}